\documentclass[11pt]{article}

\usepackage{amssymb,amsmath,latexsym}

\input psfig.sty
\psfull

\begin{document}

\newtheorem{lemma}{Lemma}
\newtheorem{theorem}{Theorem}
\newtheorem{corollary}{Corollary}
\newtheorem{definition}{Definition}
\newtheorem{example}{Example}
\newtheorem{proposition}{Proposition}
\newtheorem{condition}{Condition}

\newcommand{\bea}{\begin{eqnarray}}
\newcommand{\eea}{\end{eqnarray}}
\newcommand{\beaa}{\begin{eqnarray*}}
\newcommand{\eeaa}{\end{eqnarray*}}
\newcommand{\ben}{\begin{enumerate}}
\newcommand{\een}{\end{enumerate}}
\newcommand{\bi}{\begin{itemize}}
\newcommand{\ei}{\end{itemize}}

\newcommand{\lip}{\langle}
\newcommand{\lan}{\langle}
\newcommand{\rip}{\rangle}
\newcommand{\ran}{\rangle}
\newcommand{\uu}{\underline}
\newcommand{\oo}{\overline}
\newcommand{\til}{\tilde}

\newcommand{\La}{\Lambda}
\newcommand{\la}{\lambda}
\newcommand{\eps}{\epsilon}
\newcommand{\vph}{\varphi}
\newcommand{\al}{\alpha}
\newcommand{\bet}{\beta}
\newcommand{\gam}{\gamma}
\newcommand{\Gam}{\Gamma}
\newcommand{\kap}{\kappa}
\newcommand{\Del}{\Delta}
\newcommand{\Th}{\Theta}
\newcommand{\s}{\sigma}
\newcommand{\sig}{\sigma}
\newcommand{\Sig}{\Sigma}
\newcommand{\del}{\delta}
\newcommand{\om}{\omega}
\newcommand{\Om}{\Omega}

\newcommand{\N}{{\mathbb N}}
\newcommand{\R}{{\mathbb R}}
\newcommand{\Z}{{\mathbb Z}}
\newcommand{\bbeta}{{\text{\boldmath$\beta$}}}
\newcommand{\bdel}{\delta}
\newcommand{\bgamma}{\gamma}
\newcommand{\bnu}{\nu}
\newcommand{\bDel}{\Delta}
\newcommand{\y}{y}
\newcommand{\x}{x}
\newcommand{\X}{X}
\newcommand{\bu}{u}
\newcommand{\e}{e}
\newcommand{\m}{m}
\newcommand{\M}{M}
\newcommand{\p}{p}
\newcommand{\q}{q}
\newcommand{\bv}{v}
\newcommand{\z}{z}

\newcommand{\calA}{{\cal A}}
\newcommand{\calB}{{\cal B}}
\newcommand{\calC}{{\cal C}}
\newcommand{\calD}{{\cal D}}
\newcommand{\calF}{{\cal F}}
\newcommand{\calG}{{\cal G}}
\newcommand{\calH}{{\cal H}}
\newcommand{\calJ}{{\cal J}}
\newcommand{\calL}{{\cal L}}
\newcommand{\calM}{{\cal M}}
\newcommand{\calP}{{\cal P}}
\newcommand{\calS}{{\cal S}}
\newcommand{\calT}{{\cal T}}
\newcommand{\calU}{{\cal U}}
\newcommand{\calV}{{\cal V}}
\newcommand{\calX}{{\cal X}}
\newcommand{\calY}{{\cal Y}}

\newcommand{\proof}{\noindent {\bf Proof:\ }}
\newcommand{\proofOf}[1]{\noindent {\bf Proof of #1:\ }}
\newcommand{\remark}{\noindent {\bf Remark:\ }}
\newcommand{\remarks}{\noindent {\bf Remarks:\ }}
\newcommand{\note}{\noindent {\bf Note:\ }}
\newcommand{\esssup}{{\rm ess}\sup}
\newcommand{\essinf}{{\rm ess}\inf}
\newcommand{\pl}{\partial}
\newcommand{\noi}{\noindent}
\newcommand{\goto}{\to}
\newcommand{\ink}{\rule{.5\baselineskip}{.55\baselineskip}}
\newcommand{\qed}{\rule{.5\baselineskip}{.55\baselineskip}}

\def\ve{\varepsilon}
\def\vr{\varrho}

\newcommand{\diag}{{\rm diag}}
\newcommand{\trace}{{\rm trace}}
\newcommand{\tr}{{\rm tr}}
\newcommand{\w}{\wedge}
\newcommand{\dint}{\int\!\!\!\int}
\newcommand{\lt}{\left}
\newcommand{\rt}{\right}
\newcommand{\dist}{{\rm dist}}

\newcommand{\policy}{{u}}
\def\OBdry{\partial_o}
\def\CBdry{\partial_c}
\newcommand{\Sfrac}[2]{{{#1}\slash {#2}}}

\newcommand{\mean}[1]{\langle#1\rangle}

\date{August 14, 2003}
\title{
Explicit solution for a network control problem\\
in the large deviation regime\thanks{This
research  was supported in part by the United States--Israel Binational
Science Foundation (BSF 1999179)}}
\author{Rami Atar\footnote{Department of Electrical Engineering,
Technion -- Israel Institute of Technology,
Haifa 32000, Israel.
{\tt atar@ee.technion.ac.il}.
Research of this author also supported
in part by the fund for promotion of research at the Technion.},
Paul Dupuis\footnote{Lefschetz Center for Dynamical Systems,
Division of Applied Mathematics,
Brown University,
Providence,  R.I.\  02912.
{\tt dupuis@dam.brown.edu}.
Research of this author also supported in part by
the National
Science Foundation (NSF-DMS-0072004, NSF-ECS-9979250) and the Army Research Office 
(DAAD19-02-1-0425).}
\ and Adam Shwartz\footnote{Department of Electrical Engineering,
Technion -- Israel Institute of Technology,
Haifa 32000, Israel.
{\tt adam@ee.technion.ac.il}.
Research of this author also supported
in part by INTAS grant 265, and in part by the fund for promotion
of research at the Technion.}\\[.2in]
 }

\maketitle

\begin{abstract}
We consider optimal control of a stochastic network, where service is controlled to prevent buffer overflow.
We use a risk-sensitive escape time criterion,
which in comparison to the ordinary escape time criteria heavily penalizes exits which occur on short time intervals.
A limit as the buffer sizes tend to infinity is considered.  In~\cite{ads} we showed that, for a large class of networks,
the limit of the normalized cost agrees with the value function of a differential game.  In this game, one player controls the service discipline 
(who to serve and whether to serve),
and the other player chooses arrival and service rates in the network.
The game's value is characterized in \cite{ads} as the unique solution
to a Hamilton-Jacobi-Bellman Partial Differential Equation (PDE). 
In the current paper we apply this general theory to the important case of a network of queues in tandem.
Our main results are: (i) the construction of an explicit solution to the corresponding PDE,  and (ii) drawing out the implications for optimal risk-sensitive and robust regulation of the network.
In particular, the following general principle can be extracted.
To avoid buffer overflow there is a natural competition between two tendencies.
One may choose to serve a particular queue, since that will help prevent its own buffer from overflowing, or one may prefer to stop service, with the goal of preventing overflow of buffers further down the line.
The solution to the PDE indicates the optimal choice between these two, specifying
the parts of the state space where each queue must be
served (so as not to lose optimality), and where it can idle.
Referring to those queues which must be served
as bottlenecks, one can use the solution to the PDE to explicitly
calculate the bottleneck queues as a function of
the system's state, in terms of a simple set of equations.
\end{abstract}


\section{Introduction}

In a previous work \cite{ads}
we considered a stochastic control problem for a Markovian queueing network
with deterministic routing, where the service stations may provide service to
one or more queues, with each queue being limited by a finite buffer. The control refers to the service
discipline at the service stations, and the cost involves the time till one of the
queues first reaches its buffer limit. Such a problem can be regarded as the control of
a Markov process up to the time it exits a domain, the domain being the rectangle
associated with the buffer sizes. The cost is chosen so as to obtain a risk-sensitive
(or rare-event) control problem: one considers $E_xe^{-c\s}$ as a criterion
to be minimized, where $c>0$ is fixed and $\s$ denotes the exit time (the time when any
one of the buffers first overflows).
Such a criterion penalizes short exit times more heavily
than ordinary escape time criteria (such as $E_x\s$, a criterion to be maximized).
While the main result of \cite{ads} is the characterization of a limiting problem as the buffer sizes tend to infinity,
the current work focuses on finding explicit solutions to this limit problem and on the
interpretation thereof.

There are at least two motivations for the use of risk-sensitive criteria when designing 
policies for the control of a network.
The first is that in many communication networks performance is measured in terms of
the occurrence of rare but critical events.
Buffer overflow is a principal example of such an event.
The second motivation follows from the connection between risk-sensitive
controls and robust controls. Indeed, as discussed in \cite{dupjampet},
the optimization of a single fixed stochastic network with respect to a 
risk-sensitive cost criteria automatically produces controls with specific and 
predictable robust properties. In particular,
these controls give good performance for a family of perturbed network models 
(where the perturbation is around the design model and the size of the 
perturbation is measured by relative entropy), and with respect to a 
corresponding ordinary (i.e., not risk-sensitive) cost.

In many problems, one considers the limit of the risk-sensitive problem as a scaling parameter 
of the system converges, in the hope that the limit model is more tractable.
This is the path followed in \cite{ads}, in the asymptotic regime
where time is accelerated and buffer lengths are enlarged by a factor $n$.
The limit of the normalized cost was characterized both as the value function of
a differential game, and as the solution to a corresponding nonlinear PDE.
Both interpretations are important and useful. 
It is the interpretation as the value function of a differential game that is key in quantifying the robust aspects of the resulting control policy.
However, because the PDE gives a necessary and sufficient characterization of the (a priori unknown) value function,
it can provide qualitative information on the structure of the value function and the optimal controls.
In particularly favorable circumstances one can go even further,
and use the PDE for exact calculation and control policy synthesis.
Indeed,
if one can by any means guess a proper parametric form of the value function,
then it is sometimes possible to use the PDE to verify that this is the correct form,
and identify the unknown parameters in the representation for the value function.
The instances in the control theory literature where this has been carried out are few and far between,
especially when the state space of the system has dimension greater than two.
However, 
important information is obtained from these instances on the structure of optimal controls.

In this paper the focus is on applying the PDE obtained in \cite{ads} in the manner just described.
In particular, we treat in detail the general case of a tandem queueing network, and construct an explicit solution to the PDE.  
In a system of queues in tandem, each server offers service to exactly one queue, and therefore
the service control refers simply to whether each station should provide service or remain idle. 
Naturally, it is important to provide service to a queue in order to keep
it from reaching the buffer limit. On the other hand, if the next buffer in line is nearly full,
it might be necessary to idle the first queue in order to keep the second from overflowing.
Besides factoring in how close all buffers are to their respective limits,
one must consider the mean service rates and the likelihood of significant deviations from those mean service rates.  Is it likely that the next queue down will stall and simply stop serving for a while?
It is also possible that one will have to look even further ahead,
and consider buffers further downstream.

Although the optimal control problem for the Markovian queueing system is fully described
by a dynamic programming equation \cite{ads}, such equations are typically solved numerically.
It is hard to extract any global structural information from the exact equation,
and even a numerical solution may not be feasible when buffers are large or the dimension of the state space is moderate. 
As discussed previously, the solution to the PDE for the limiting problem turns out to simplify things
significantly. Roughly speaking, the PDE indicates the following structure of
the asymptotic optimally controlled
network. In the interior of the domain
(i.e., when all buffers are away from their maximum capacity), and depending on the state of the system,
service must be provided at
certain service stations, while other stations may idle without causing loss of optimality.
In that sense, the limit problem sharpens the control policy by emphasizing
the importance of serving those `bottleneck' buffers
(it is crucial to serve the bottleneck buffers, and completely
unimportant to serve the others).
The identification of the bottlenecks is nontrivial, and indeed queues with the smallest service rate are not necessarily  bottlenecks.
Instead, as hinted above, identification of the bottlenecks must include consideration of at least the following: (i) the relative closeness of all buffers to their maximum value,
(ii) relative mean service rates,
and (iii) relative uncertainties in the service rates.
We identify a system of algebraic equations whose solution identifies all bottlenecks.

There is relatively little work on risk-sensitive and robust control of networks.
Ball et.\ al.\ have considered a robust formulation for network problems
arising in vehicular traffic \cite{balday2}, where the cost structure is qualitatively
different. Dupuis studies a robust control problem for networks in a deterministic setting
and obtains explicit solutions for the value function \cite{dup}.
The cost there is, in a sense, antipodal to the one considered in the current work,
namely the time till the system becomes empty (a criterion to be minimized).
For other recent work on queueing control in a large deviation regime
see Stolyar and Ramanan \cite{storam}, where a
single server non-Markovian system (with quite general stationary increments input flows)
is studied, and a particularly simple scheduling control policy is shown
to be asymptotically optimal (see also an extension of the work in
Stolyar~\cite{sto}).

The paper is organized as follows.
Section 2 introduces the model and the PDE and states the main result.
Section 3 contains discussion and interpretation.
In Section 4 we prove the main result.

Notation:
The symbol $\vee $ stands for
maximum, while $\w $ stands for minimum. Denote scalar product between two
vectors as $ \x \cdot \y $.
For integers $i\le j$, let $[i,j]\doteq\{1,\ldots,j\}$.

\section{Model and preliminaries}

\noi\uu{\bf The queueing network control problem.}
We consider the following tandem network.
There are $J$ queues and $J$ servers. Customers present at queue $i$
at a certain time are said to be of class $i$ at that time.
Customers arrive to queue $1$ according to a
Poisson process of rate $\la\ge0$. Service at queue $i$ is provided by
server $i$ at exponential time
with parameter $ \mu_i>0$, mutually independent and independent of the
arrivals.
After a customer is served by server $i$, it moves to queue
$r(i)$, where $r(i) = i+1 $, $i=1,\ldots,J-1$, $ r(J) = 0$,
and $i=0$ is used to denote the ``outside''.
The state of the network is the vector of queue lengths, denoted
by $\X$.
Let $\{ \e_i;i=1,\ldots,J\} $ denote the unit coordinate vectors,
let $\e_0 = 0 $ and denote
\begin{equation}\label{egam}
\bgamma_i=\e_i-\e_{r(i)}.
\end{equation}
Note that $\bgamma_i=\e_i-\e_{i+1}$ for $i=1,\ldots,J-1$, $\bgamma_J
=\e_J$, and that
following service to queue $i$ the state changes by $-\bgamma_i$.
The control is specified by the vector $ u =(u_1,...,u_J)$,
where $u_i=1$ if customers in queue $i$ are given service and $u_i=0$
otherwise.
We next consider the scaled process $\X^n$ under the scaling
which accelerates time by a factor of $n$ and down-scales space by
the same factor. We are interested in a
risk-sensitive cost functional that is associated with
exit from a bounded set.
Let $G$ be the rectangle defined through
\begin{equation}\label{def:G}
 G=\{(x_1,\ldots,x_J): 0\le x_1<z_1;\ 0\le x_i\le z_i, \,
i = 2 , \ldots , J \},
\end{equation}
for some $z_i>0$, $i=1,\ldots,J$.
Note that $G$ contains parts of, but not all of its boundary.
Let
\[
\sigma^n \doteq \inf \{t:\X^n(t) \not \in G \}.
\]
The control problem is to minimize the cost
$E_\x e^{-nc\sigma^n}$, where $E_\x$ denotes expectation starting
from $\x$, and $c>0$ is a constant.
With this cost structure ``risk-sensitivity'' means that atypically 
short exit times are weighted heavily by the cost.
As a result, even if exit within short time
occurs with small probability, it may have a noticeable effect on the cost,
and thus a ``good'' control will attempt to avoid such events
as much as possible.

Although this paper do not directly treat the
stochastic control problem but the corresponding PDE, we give here the precise
formulation of the former, for completeness.
For each $n$, the state space is $G^n \doteq n^{-1}\Z_+^J \cap G$. 
The control (or action) space is given by 
$$
 U \doteq \left\{(u_i), i=1,\ldots,J: 0 \le u_i\le1, \
 i=1,\ldots,J\right\}.
$$
One considers state processes $\X^n$ that are similar
to those defined for the original queueing network,
except that time is accelerated (equivalently, arrival and service rates
are multiplied) by a scaling parameter $n$, and space is scaled
down by the same factor. More precisely, for $u\in U$ and $f:\Z_+^J\to\R$,
let
$$
\calL^uf(\x)=\la [f(\x+\e_1)-f(\x)]
+\sum_{i=1}^Ju_i\mu_i1_{\{\x-\bgamma_i\in\Z_+^J\}}[f(\x-\bgamma_i)-f(\x)],
\quad \x\in\Z_+^J.
$$
For $n\in\N$ let
$$
\calL^{n,u}f(\x)=n\calL^ug(n\x),
$$
where $f:n^{-1}\Z_+^J\to\R$ and $g(\cdot)=f(n^{-1}\cdot)$.
A {\em controlled Markov process} starting from $\x\in G^n$
will consist of a complete filtered
probability space $(\Om,\calF,(\calF_t),P_x^{n,u})$,
a state process $\X^n$ taking values in  $G^n$ and
a control process $u$ taking values in  $U$,
such that $\X^n$ is adapted to $\calF_t$,
$u$ is measurable and adapted to $\calF_t$,
$P_\x^{u,n}(\X^n(0)=\x)=1$,
and for every function $f:G^n\to\R$
\[
f(\X^n(t))-\int_0^t\calL^{n,u(s)}f(\X^n(s))ds
\]
is an $\calF_t$-martingale.
$E_\x^{n,u}$ denotes expectation with respect to $P_\x^{n,u}$.
For a parameter $c>0$,
the value function for the stochastic control problem is 
defined by
\begin{equation}\label{eq:control}
  V^n(\x) \doteq -\inf n^{-1}\log E_\x^{u,n}e^{-nc\sig_n},
  \quad \x\in G^n,
\end{equation}
where the infimum is over all controlled Markov processes.

\begin{figure}
\centerline{
\begin{tabular}{cc}
\psfig{file=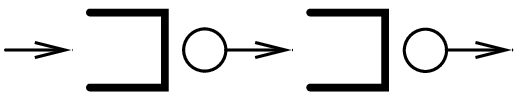}\qquad\qquad\qquad
\psfig{file=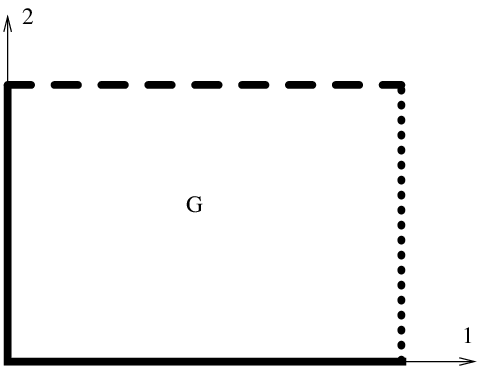}
\end{tabular}
}
\caption{
A simple queueing network, a rectangular domain and three
types of boundary.
Full line: $~$ \hfill $\pl_+G$, dashed line: $\pl_cG$,
and dotted line: $\pl_oG$
\hfill $~$ }
\label{fig1}
\end{figure}

\noi\uu{\bf The domain and its boundary.}
It is possible for the controller to prevent any but the first queue
from exceeding $ z_i $, simply by turning off service to the preceding
queue. However, the controller cannot prevent overflow of the first queue.
Although it is in principle possible that the dynamics
could exit through the portion of the boundary defined by queues
$ 2 , \ldots , J$, it is always optimal for
the controller to not allow this.
Consider the simple two-class network illustrated in Figure \ref{fig1}.
The controller can prevent exit through the dashed portion of the 
boundary simply by stopping service at the first queue. 
As a consequence, there are in general three different types
of boundary--the constraining boundary due to non-negativity constraints on queue length,
the part of the boundary where exit can be blocked, and the remainder.
These three types of boundary behavior result
in the PDE in three types of boundary conditions.
The three portions of the boundary are explicitly given as
$$
\pl_c G=\{(x_1,\ldots,x_J): 0\le x_1<z_1,\ \mbox{ and }\ x_i=z_i\ \mbox{for
some}\,i> 1\}
$$
$$
\pl_o G=\{(x_1,\ldots,x_J): x_1=z_1, \ \mbox{ and}\ 0\le x_i\le z_i\ \mbox{for
all}\ i>1\}
$$
$$
\pl_+G= \{(x_1,\ldots,x_J):
\mbox{$x_i<z_i$ for all $i$, and $x_i=0$ for some $i$}\}.
$$
Note that $\pl_cG$, $\pl_oG$ and $\pl_+G$
partition the boundary $\pl G$ of $G$.
Also, $\pl_cG\subset G$ while $\pl_oG\cap G=\emptyset$.
As usual, we will denote $G^o=G\setminus\pl G$
and $\bar G=G\cup\pl G$.

\noi\uu{\bf The Hamiltonian, PDE and viscosity solutions.}
It is a standard fact that the value functions $V^n$
of the stochastic control problem considered above can be characterized
by a dynamic programming equation. The results of \cite{ads}
show that also the limit of $V^n$ as $n\to\infty$
can be characterized by a certain equation, namely a Hamilton-Jacobi-Bellman
PDE, that can in fact be regarded as the dynamic programming equation
for a certain deterministic game problem. For details on this PDE
we need some notation (the reader is referred to \cite{ads} for details
on the game).
Let $l:\R\to\R_+\cup\{+\infty\}$
be defined as
$$
 l(x)\doteq \lt\{\begin{array}{ll}x\log x-x+1 & x\ge0,  \\ +\infty & x<0,
     \end{array}\rt.
$$
where $0\log0\doteq0$.
Define
\begin{equation}
\label{em}
 \M =\{m=(\bar\la,\bar\mu_1,\ldots,\bar\mu_J):
\bar\la \geq 0 ,\ \bar \mu_i \geq 0 \}.
\end{equation}
For $ u\in U$ and $ m\in \M$ define
\begin{align*}
v ( u, m)    & \doteq \bar \lambda \e_1 - \sum_{i=1}^J
               u_i \bar\mu_i \bgamma_i, \\
\rho( u, m) & \doteq \lambda l\left( \frac{\bar{\lambda}}{
\lambda}\right) +\sum_{i=1}^{J}u_i\mu _il\left( \frac{\bar{\mu}_i}
{\mu_i}\right),
\end{align*}
where $\bgamma_i$ are as in \eqref{egam}.
Let
\begin{align}\label{def:H1}
H(p, u, m) &= c+ p\cdot v ( u, m)+\rho( u, m)\\ \notag
&=
c+
\lt[\bar\la p_1+\la l\lt(\frac{\bar\la}{\la}\rt)\rt]
+\sum_{i=1}^J u_i
\lt[-\bar\mu_i\bgamma_i\cdot \p+\mu_i l\lt(\frac{\bar\mu_i}{\mu_i}\rt)\rt]
\end{align}
and let the Hamiltonian be defined as
\begin{equation}\label{def:H}
 H(p)=\sup_{ u\in U}\inf_{ m\in\M}H(p,u,m)=\inf_{m\in \M}\sup_{u\in U}
H(p, u, m),
\end{equation}
where the last identity, expressing the {\it Isaacs Condition} is proved in \cite{ads}.
The following simplification in the structure of 
the Hamiltonian will be useful.
Using convexity and the fact that the slope of $l$ at $0^+$ is
$-\infty$, the minimum over $ m$ is attained
at $\bar\la=\la e^{-p_1}$, $\bar\mu_i=\mu_i e^{\bgamma_i\cdot \p}$.
A straightforward calculation then shows that
\begin{align}\notag
H(\p, u)\doteq \inf_m H(\p, u, m) &=
c+[(\la-\bar\la)+u_i(\mu_i-\bar\mu_i)]\\
&=
c+\la(1-e^{-p_1})
+\sum_{i=1}^Ju_i\mu_i(1-e^{\bgamma_i\cdot \p}).
\label{eq:116}
\end{align}
Define $ I (\x) = \{ i : x_i = 0 \} $.
The PDE of interest is the following.
\begin{equation}\label{eq:pde}
\lt\{
\begin{array}{ll}
 H(DV(\x))=0, & \x\in G^o,\\ \\
 DV(\x)\cdot \bgamma_i=0, & i\in I(\x),\ \x\in\pl_+G ,\\ \\
 V(\x)=0, & \x\in\OBdry G.
\end{array}\rt.
\end{equation}
Since typically such equations do not possess classical solutions,
the framework of viscosity solutions is useful (see \cite{barcap}).
This framework allows
for functions that are merely continuous to be regarded as solutions.
For $\x\in G$,
the set of superdifferentials $D^+V(\x)$ and the set of subdifferentials
$D^-V(\x)$ are defined as
\begin{equation}\label{def:Dplus}
D^+V(\x)=\lt\{\p:\limsup_{\y\to \x}\frac{V(\y)-V(\x)-\p\cdot(\y-\x)}{|\y-\x|}\le0\rt\},
\end{equation}
\begin{equation}\label{def:Dminus}
D^-V(\x)=\lt\{\p:\liminf_{\y\to \x}\frac{V(\y)-V(\x)-\p\cdot(\y-\x)}{|\y-\x|}\ge0\rt\}.
\end{equation}

\begin{definition}
$V$ is a viscosity solution to (\ref{eq:pde}) if
\begin{equation}\label{eq:up}
H(\p)\vee\max_{i\in I(\x)}\p\cdot\bgamma_i\ge0, \quad
\p\in D^+V(\x),\ \x\in G,
\end{equation}
\begin{equation}\label{eq:down}
H(\p)\w\min_{i\in I(\x)}\p\cdot\bgamma_i\le0, \quad
\p\in D^-V(\x),\ \x\in G\setminus\pl_cG,
\end{equation}
and $V(\x)=0$ for $\x\in\pl_oG$.
\end{definition}
The following is a special case of Theorem 2 of \cite{ads}.
\begin{theorem}\label{th:main}
There exists a unique Lipschitz viscosity solution $V$ to (\ref{eq:pde}).
Moreover,
if $\x_n\in G^n$, $n\in\N$ are such that $\x_n\to \x\in G$, then
$\lim_{n\to\infty}V^n(\x_n) = V(\x)$.
\end{theorem}
For $i=1,\ldots,J$, let $\beta_i$ denote the unique positive solution to
\begin{equation}\label{eq:100}
c+\la(1-e^{\beta_i})+\mu_i(1-e^{-\beta_i})=0.
\end{equation}
Set $b_i=\beta_i\sum_{j=1}^i\e_j$.
Our main result is the following.
\begin{theorem}\label{th:tandem}
Assume $c>0$.
Then the viscosity solution to the PDE (\ref{eq:pde}) is given by
\begin{equation}\label{eq:101}
V(\x )=\min_{i=1,\ldots,J}b_i\cdot(\z -\x ).
\end{equation}
\end{theorem}

\section{Some remarks}\label{secrems}

\subsection{Interpretation}

As commented above, it is possible to characterize
the value function $V^n$
as the solution to a dynamic programming equation on the discrete space
$G^n$. Moreover, if $V^n$ is available in explicit form, one can
use the dynamic programming equation to specify an optimal control
policy for the problem. The results above provide an explicit
expression only to the limit $V=\lim_n V^n$, and therefore we are unable
to specify an optimal policy for the problem with finite $n$.
Instead, we shall use the quantity
$V$ along with the HJB equation \eqref{eq:pde} to propose a policy for
the queueing network problem (and the corresponding value $V^n$).
It is plausible that the proposed policy is, in a sense,
asymptotically optimal as $n\to\infty$ but we do not attempt
to prove such a statement here.
Instead, the discussion below could be regarded as a natural
interpretation of the result.
Moreover, since $V$ can be considered as a good approximation to
$V^n$ (for large $n$), it can be used as an initial condition
in a value iteration procedure that calculates $V^n$ and an optimal
policy.

The solution to the PDE stated here may be interpreted as follows.
Let $\x $ be an interior point of $G$, and assume that
$V$ is differentiable at $x$.
The form \eqref{eq:101} implies that for some $j=j(\x)$, 
\begin{equation}\label{e300}
DV(\x)=-b_j
=-\beta_j\sum_{i=1}^j\e_i \, .
\end{equation}
Hence by \eqref{eq:116}, the equation $H(-b_j)=0$ has the form
$$
\sup_{u\in U}\lt[c+\la(1-e^{\beta_j})+u_j\mu_j(1-e^{-\beta_j})\rt]=0.
$$
It is seen that the supremum is attained at any $u$ for which $u_j=1$
and $u_i\in[0,1]$, $i\ne j$. The interpretation in terms of the service policy
is that it is important to serve class $j$, while optimality does not depend
on the service given to other classes.
{\em The station $j=j(x)$ can therefore be regarded as a bottleneck:}
it is crucial to serve station $j$ when at state $x$.
Since the PDE describes the limit of the stochastic control problem,
one expects
the bottleneck stations to have a similar property in the stochastic problem:
Although it may be optimal for all servers to not idle when the system is near
an interior point $\x$, it is significant for the bottleneck server
to work while if the other servers idle the cost is affected only by little.

A closer look at \eqref{e300} reveals that all bottleneck stations
belong to a certain set $A'$ defined below.
More precisely,
given an interior point $x$ where $V$ is differentiable
and denoting $y=z-x$, using \eqref{eq:101}, we see
that \eqref{e300} holds for $j$ if and only if $j$ satisfies
\begin{equation}\label{e301}
\beta_j(y_1+\cdots+y_j)\le\beta_i(y_1+\cdots+y_i)\quad
\forall\ i\ne j.
\end{equation}
We claim that
a necessary condition for $j$ to satisfy \eqref{e301} (for some $y$)
is $j\in A'$, where
$$
A'=\{k\in[1,J]:\mu_k\le\mu_l,\ \mbox{ for all } l<k\},
$$
and by convention, $1\in A'$.
To this end, write the explicit form of the solution to~\eqref{eq:100} as
\begin{equation*}
e^{\beta_i } = (2\la)^{-1}[ c + \lambda + \mu_i
                   + ( c^2 + \lambda^2 + \mu_i^2 
                       + 2 c ( \lambda + \mu_i) )^{1/2}] \, ,
\end{equation*}
and note that the positive solutions $\beta_i$ to
\eqref{eq:100} are monotone in $\mu_i$, in the sense that
\begin{equation}\label{eq:mon}
\mu_i<\mu_j\quad \Longleftrightarrow\quad \beta_i<\beta_j.
\end{equation}
If \eqref{e301} holds then $\beta_i\le\beta_j$ for $j=1,\ldots,i$,
and it follows from \eqref{eq:mon} that $j\in A'$.

It should be emphasized that the condition $j\in A'$ is only
necessary for $j$ to be a bottleneck, and being a bottleneck is
really a function of the system's state. A more complicated necessary
condition is given in Section \ref{sec:proof}, which is of the form
$j\in A(x)$.
One interesting and perhaps counterintuitive phenomenon is that it is possible
that a station is nearly empty while others are far from being empty,
and still it is a bottleneck.

\begin{figure}
\centerline{
\begin{tabular}{cc}
\psfig{file=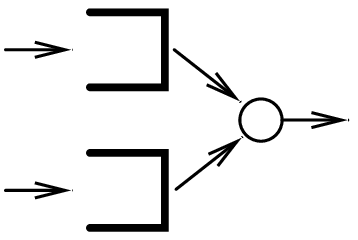}\qquad\qquad\qquad
\psfig{file=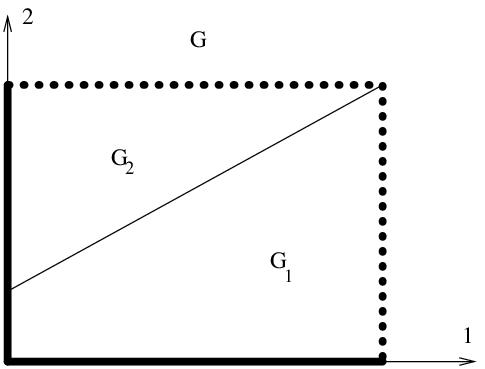}
\end{tabular}
}
\caption{Priority to class $i$ when the state is in $G_i$, $i=1,2$.}
\label{fig:example}
\end{figure}

In the rest of this subsection we discuss a different queueing system
for which a similar analysis is possible, including an explicit form
for the limit of the value function and a clear interpretation of
it regarding nearly optimal service policies for the queueing system
(see \cite{ads} for proof of the results quoted below).
Consider a single server that provides
service to $J$ classes (each customer requires service once).
Service rate to queue $i$ is exponential with rate $ \mu_i $,
and the arrival process to class $i$ is Poisson of parameter $\la_i>0$,
$i=1,\ldots,J$. All arrival and service processes are mutually independent.
One considers the control space
\[
U=\left \{ u: \sum_{i=1}^J u_i\le 1, u_i \ge 0, i=1,\ldots,J \right \},
\]
where $u_i$ represents the fraction of service allocated to class $i$
by the single server.
Analogously to the problem discussed above, one defines a domain
$G=\{(x_1,\ldots,x_J):0\le x_i< z_i,1=1,\ldots,J\}$, scaled queueing
processes $X^n$ and exit time $\s^n$, and attempts minimizing the cost
$E_xe^{-nc\s^n}$ over an appropriate class of control processes.
Let $V^n(x)$ denote the infimum.
One can then show that $\lim_nV^n(x_n)$
exits whenever $x_n\to x\in G$.
Moreover, the limit $V$ is characterized in terms of a PDE
of the form \eqref{eq:pde}, where now $\gamma_i=e_i$, $i=1,\ldots,J$,
and $H(p)=\sup_{u\in U}H(p,u)$. Here $H(p,u)$ has the form
(analogous to \eqref{eq:116}):
\begin{equation}\label{216}
H(\p,u)=
c+\sum_{i=1}^J[\la_i(1-e^{-p_i})+u_i\mu_i(1-e^{p_i})],
\end{equation}
and the limit $V(x)$
can be explicitly calculated provided $c$ is large enough
as follows.
For $i=1,\ldots,J$, define $ \alpha_i $ as the unique positive
solution to $c+\la_i(1-e^{\al_i})+\mu_i(1-e^{-\al_i})=0$, namely
$$
e^{\alpha_i } = (2\la)^{-1}[ c + \lambda + \mu_i
                   + (( c + \lambda + \mu_i )^2
		   - 4 \lambda \mu_i)^{1/2}]\, .
$$
Then if $c$ is large enough one has
\begin{equation}\label{eq:form}
V(\x )=\min_{i=1,\ldots,J}\al_i(z_i-x_i).
\end{equation}
Let $\x$ be a point in the interior of $G$
where the gradient of $V$ is well defined.
Then the PDE is satisfied in the classical sense at this point.
Thus, $H(DV(\x))=H(-\al_j\e_j)=0$, where $j=j(\x)$.
Using \eqref{216}, the equation $H(-\al_i\e_j)=0$ takes the form
$$
\sup_{u\in U}\lt[c+\la_j(1-e^{\al_j})+
u_j\mu_j(1-e^{-\al_j})\rt]=0.
$$
Clearly the supremum is attained at $u_i=1_{i=j}$,
$i=1,\ldots,J$.
This means that it is optimal to serve class $j(\x)$ at the
state $\x$.
In the totally symmetric case, where $\mu_i=\mu$, $\la_i=\la$,
$z_i=z$ for all $i$, the solution takes the form
$V(\x )=\al\min_i (z-x_i)$, and the optimal
service discipline can be interpreted as ``serve the longest
queue.'' An asymmetric two dimensional example
is given in Figure \ref{fig:example}, where the domain $G$
is divided into two subdomains $G_1$ and $G_2$ in accordance
with the structure (\ref{eq:form}), and the optimal service
discipline corresponds to giving priority to
class $i$ when the state is within $G_i$, $i=1,2$.
This discipline gives priority {\em to the queue with the (weighted)
shortest free buffer space.}

\subsection{The perturbed rates}

In the asymptotic analysis of rare events it is often the case
that most of the probability mass of an event of interest
is concentrated on the event that the stochastic process, say $X^n$,
nearly follows a certain deterministic path, for large values of $n$.
For example,
the most likely way a stable M/M/1 queue overflows is by nearly
following a deterministic path that is the solution to a fluid model
in which arrival and service rates are reversed.
This  appears as a consequence of time-reversal arguments in~\cite{sw}.
For many other cases where one can compute
deterministic paths about which most mass is concentrated see \cite{swbook}.
Similarly, in the current stochastic control problem,
when the system operates under the optimal control most
contribution to the cost is obtained when the process $X^n$ nearly
follows a certain deterministic path. This path is now
the solution to a fluid model with a control $u$ and perturbed
arrival and service rates $m=(\bar\la,\bar\mu_1,\ldots,\bar\mu_J)$.
At a point $x$ where the viscosity solution $V$ to \eqref{eq:pde}
is differentiable, the correct values of $u$ and $m$ are those that achieve
the max-min in $H(DV(x))=\sup_u\inf_mH(DV(x),u,m)$ (cf.\ \eqref{def:H}).

We remark that in the example of
tandem queues the optimal perturbed rates
satisfy certain relations with the unperturbed rates, as shown in
the two equations below.
By \eqref{eq:116}, the relation $H(DV)=0$ implies
\begin{equation}\label{e302}
\bar\la+\sum_{i=1}^Ju_i\bar\mu_i=c+\la+\sum_{i=1}^Ju_i\mu_i \, .
\end{equation}
Moreover, by \eqref{eq:116}, by the fact that the minimum over $m$
is attained at
$\bar\la=\la e^{-p_1}$, $\bar\mu_i=\mu_i e^{\bgamma_i\cdot \p}$,
and using $\sum_{i=1}^J \bgamma_i = \e_1$, we have
\begin{equation}\label{e303}
\bar\la\prod_{i=1}^J\bar\mu_i=\la\prod_{i=1}^J\mu_i.
\end{equation}
Relation \eqref{e302}
was noticed by Avram \cite{avr}.
The relation \eqref{e303} appears to be new.
A generalization of these equations to a more general network
is possible and will appear elsewhere.

\remark An interesting relation between the roots $\al_i$ and $\beta_i$
and the busy cycle period was
pointed out to us by Boxma \cite{box}.
For example, if $B$ denotes the busy cycle period for an M/M/1
queue under the stationary distribution, then $Ee^{cB}=e^\beta$,
where
$$
c+\la(1-e^{\beta})+\mu(1-e^{-\beta})=0
$$
(compare with \eqref{eq:100}).

\section{Proof of the main result}\label{sec:proof}

Recall that $I(\x )=\{i\in[1,J]:x_i=0\}$, and denote
$B(\x )=\{i\in[1,J]:x_i=z_i\}$. Note that $I( \x )$ and $B(\x )$
do not intersect. Since points $\x $ for which $x_1=z_1$ are not in $G$,
we have $1\not\in B(\x )$ for all $\x \in G$.

We show that the minimum over $i=1,\ldots,J$
in \eqref{eq:101} can equivalently be performed over a smaller
(state-dependent) set. 
As in Section \ref{secrems}, let 
\begin{equation}\label{eq:112}
A'=\{i\in[1,J]:\mu_i\le\mu_j,\ \mbox{ for all } j<i\},
\end{equation}
and by convention let $1\in A'$.
Let also
$$
A(x)=\{i\in A':\text{ either $i+1\not\in B(\x)$; or 
$j>i$ and $[i+1,j]\subset B(\x )$
imply $\mu_i<\mu_j$}\}.
$$

\begin{lemma}\label{lem:1}
The minimum in~\eqref{eq:101} is obtained over the indices in the
set $A(\x)$ defined above.
More precisely,
\[
V(\x)=\min_{i=1,\ldots,J}b_i\cdot(\z-\x)
=\min_{i\in A(\x)}b_i\cdot(\z-\x) \, .
\]
\end{lemma}
\proof
Fix $ \x \in G $ and recall that $ z_1 - x_1 > 0 $.
If
\begin{equation}\label{eq:mon1}
b_i\cdot(\z -\x )\le\min_{j\ne i}
b_j\cdot(\z -\x ),
\end{equation}
then $\beta_i\le\beta_j$ for $j \in [1,i]$,
and it follows from \eqref{eq:mon} that $i\in A'$.
Thus
\[
V (\x ) = \min_{i=1,\ldots,J}b_i\cdot(\z -\x )
=\min_{i\in A'}b_i\cdot(\z -\x )
\]
for all $\x \in G $ (although not every $ i \in A' $ is a minimizer).
Fix now $i \in A'$, and suppose that
\begin{equation}\label{eq:assume}
\mbox{there is some $j>i$ so that 
$[i+1,j]\subset B(\x )$ and $\mu_i \geq \mu_j$.}
\end{equation}
Since $ i \in A'$, if $ i+1 \not\in A' $ then $ \mu_{i+1} > \mu_i $.
But then,
if in addition $ i+2 \not\in A' $ we have $ \mu_{i+2} > \mu_i $.
So, if $ k \not\in A'$ for all $ i+1 \le k \le j-1 $ then
$ \mu_k > \mu_i $, and together with $ \mu_j \le \mu_i $ we have
that $ j \in A' $. We conclude that under~\eqref{eq:assume},
\begin{equation}\label{eq:200}
\mbox{there is $k\in[i+1,j]$ such that $k\in A'$.}
\end{equation}
By~\eqref{eq:assume},
$\beta_i \geq \beta_j$ and so
\begin{equation}\label{eq:201}
b_i \cdot(\z -\x )\geq b_j \cdot(\z -\x ).
\end{equation}
Therefore, for $k$ as in \eqref{eq:200} and under the assumptions
in~\eqref{eq:assume},
\[
b_i\cdot(\z-\x)\ge b_k\cdot(\z-\x) ,
\]
and since $k\in A'$, we need not consider $i$ in the minimum.
As a result,
\eqref{eq:101} is equivalently given as
$$
V(\x )=\min_{i\in A(\x )}b_i\cdot(\z -\x ),\quad \x \in G,
$$
where $ A (\x ) $ is obtained from $ A' (\x ) $ by deleting those
$ i \in A' (\x ) $ which are followed by empty queues,
$ [ i+1 , j ] \subset B (\x ) $, and such that $\mu_i \ge \mu_j $.
This leaves in $ A (\x ) $ an index $i$ from $ A' (\x ) $
only if either $ z_{i+1} > x_{i+1} $, or 
$ z_{k} = x_{k} $ for $ i+1 < k \le j $ implies $ \mu_i < \mu_j $.
More formally,
\begin{align}
\notag
A(\x) &=
A'\setminus\{i\in A':\exists\ j>i,\ [i+1,j]\subset B(\x),\ \mu_i\ge\mu_j\}\\
\label{eq:113}
&=
\{i\in A':\text{ either $i+1\not\in B(\x)$; or 
$j>i$ and $[i+1,j]\subset B(\x )$
imply $\mu_i<\mu_j$}\}.
\end{align}
\qed

The form of the proposed solution is the minimum of smooth functions.
The set of superdifferentials
at a point where the function is not smooth
can be seen (using the definition \eqref{def:Dplus})
to consist of the convex hull of the gradients of the
smooth functions defining it, at that point. Also,
at the boundary, the fact that there are less constraints introduced
by \eqref{def:Dplus} on $p\in D^+V$ than there are when $x$
is in the interior, has an effect of enlarging the set further.
Thus, for example, the set of
superdifferentials of the zero function from $\R_+$ to $\R$ at zero
is $\R_+$.
Using these considerations and the rectangular structure
of the domain, one finds the general form for the superdifferential
as follows.
Let $\x \in G$ be fixed. Set $A=A(\x )$, $I=I(\x )$, $B=B(\x )$, and
$O=\{1,\ldots,J\}\setminus (I\cup B)$.
Then any element $\p \in D^+V(\x )$ is given as
$$
\p =-\sum_{i\in A}\nu_ib_i+\bdel,
$$
for some $ \{\nu_i \}$, where
\begin{equation}\label{eq:102}
\del_i\ge0,\ i\in I,\quad
\del_i\le0,\ i\in B,\quad
\del_i=0,\ i\in O,
\end{equation}
and
\begin{equation}\label{eq:103}
\nu_i\ge0,\ i\in A,\quad \nu_i=0,\ i\not\in A,\quad \sum_{i\in A}\nu_i=1.
\end{equation}
If we denote
\begin{equation}\label{eq:105}
\del_{J+1}=0,
\end{equation}
then, using \eqref{eq:116}, the Hamiltonian is expressed as
\begin{align}\label{eq:114}
\notag
H(\p ) &= c+\la(1-e^{-p_1})+\sum_{i=1}^J0\vee\mu_i(1-e^{\p \cdot\bgamma_i})\\
     &= c+\la(1-e^{\sum_{i  \in A}\nu_i\beta_i-\del_1})
     +\sum_{i=1}^J0\vee\mu_i(1-e^{-\nu_i\beta_i+\del_i-\del_{i+1}}).
\end{align}

\proofOf{Theorem~\ref{th:tandem}}

\noi\uu{\it Verifying the PDE for superdifferentials.}
To verify \eqref{eq:up}, it suffices to show that
$H(\p )\ge0$ whenever
$\p \in DV^+(\x )$ and 
$\p \cdot\bgamma_i<0$ for all $i\in I$.

\noi\uu{\it Step 1.} We show first that $H(\p )\ge0$ whenever
$\p \in DV^+(\x )$ and  $\p \cdot\bgamma_i\le0$ for all $i=1,\ldots,J$.
The forms of $b_i$ and $\bgamma_i=\e _i-\e_{i+1}$ imply the last inequality can be rewritten
\begin{equation}\label{eq:104}
-\nu_i\beta_i+\del_i-\del_{i+1}\le0,\ i=1,\ldots,J.
\end{equation}
In this case, \eqref{eq:114} becomes
\begin{equation}\label{eq:109}
H(\p ) = h(\bnu,\bdel)
 \doteq c+\la(1-e^{\sum_{i \in A}\nu_i\beta_i-\del_1})
+\sum_{i=1}^J\mu_i(1-e^{-\nu_i\beta_i+\del_i-\del_{i+1}}).
\end{equation}
The constraints we have put on $(\bnu,\bdel)$ define a convex
set $S$, namely, the set of $(\bnu,\bdel)\in\R^J\times\R^{J+1}$ satisfying
\eqref{eq:102}, \eqref{eq:103}, \eqref{eq:105} and \eqref{eq:104}.
The set $S$ is bounded for the following reasons.
First, since $\bnu$ is a probability vector (see \eqref{eq:103}) we have $\nu_i\in[0,1]$.
Next, by \eqref{eq:102} and $1\not\in B$,
 $\del_1\ge0$.
Finally,
it follows from \eqref{eq:104}
(using the convention $\sum_s^t=0$ if $t<s$) that 
$$
\del_1-\sum_{j=1}^{i-1}\nu_j\beta_j
\le\del_i\le\sum_{j=i}^J\nu_j\beta_j,
\quad i=1,\ldots, J.
$$
This shows the boundedness.

The function $(\bnu,\bdel)\mapsto h(\bnu,\bdel)$ is concave, and therefore to prove $h\ge0$ on $S$ it suffices to check the inequality
on the set $E$ of extreme points of $S$, which contains a finite number of
points since the constraints are linear.
\begin{lemma}
All points $(\bnu,\bdel)\in E$ have the following form:
\begin{align*}
\nu_k & = 1 \quad \text{for some $k$, and} \\
\del_i^{(r,k)} & =\beta_k(1_{i\ge r+1}-1_{i\ge k+1}),\quad
i=1,\ldots,J,\ r=s,\ldots,t-1,
\end{align*}
where $ t, s $ are defined below.
\end{lemma}
\proof
We will first obtain the general form of $\bnu$
and then that of $\bdel$.

We claim that for any $(\bnu,\bdel)\in E$, $\bnu$ is of the form
$\bnu=1_k$ (short for $\nu_i=1_{i=k}$, $i=1,\ldots,J$),
for some $k\in A$. Assume this is false.
Then there are $l,m\in A$, $l<m$, for which $\nu_l,\nu_m\in(0,1)$.
We will show that
there is a vector $\bDel_m$, such that replacing $(\bnu,\bdel)$ by
$(\bnu',\bdel')=(\bnu-\eps \e_m,\bdel+\eps\bDel_m)$ maintains
the relations \eqref{eq:102}, \eqref{eq:105} and \eqref{eq:104}
for both $\eps>0$ and $\eps<0$ (provided $|\eps|$ is small).
A similar statement will hold also for
$(\bnu'',\bdel'')=(\bnu+\eps \e_l-\eps \e_m,\bdel-\eps\bDel_l+\eps\bDel_m)$,
and as a result all of \eqref{eq:102}, \eqref{eq:103},
\eqref{eq:105} and \eqref{eq:104}
will hold for both $\eps>0$ and $\eps<0$, a contradiction to
$(\bnu,\bdel)$ being an extreme point of $S$.
The construction of $\bDel_m$ based on $(\bnu,\bdel)$
can be mimicked to construct $\bDel_l$ based on $(\bnu',\bdel')$
(in particular, no use is made of the fact that $\nu_i$ sum to one,
but only that some of its components are within $(0,1)$),
and therefore the latter construction is omitted.

Case 1: {\it Inequality \eqref{eq:104} holds as a strict inequality
for $i=m$, i.e., $-\nu_m\beta_m+\del_m-\del_{m+1}<0$.}
Then $-(\nu_m-\eps)\beta_m+\del_m-\del_{m+1}\le0$
also holds (for both $\eps$ positive and negative),
provided that $|\eps|$ is small. Here we take $\bDel_m=0$.
We see that $(\bnu-\eps \e_m,\bdel)$ satisfies the requirements.

Case 2: {\it Inequality \eqref{eq:104} holds with equality
for $i=m$, i.e., $-\nu_m\beta_m+\del_m-\del_{m+1}=0$.}
Since $\nu_m\beta_m>0$, either $\del_m>0$ or $\del_{m+1}<0$.
Assume $\del_m>0$ (the case $\del_{m+1}<0$ can be treated analogously,
and is therefore omitted).
Let $m'$ be the smallest $j\in[1,m]$ for which \eqref{eq:104}
holds with equality for all $i\in[j,m]$. For $i\in[m',m-1]$
(the set being empty and the statement void if $m'=m$)
we have $-\nu_i\beta_i+\del_i-\del_{i+1}=0$, hence if
$\del_{i+1}>0$ then $\del_i>0$. Together with the fact that
$\del_m>0$, this shows that
\begin{equation}\label{eq:108}
\del_i>0,\ i\in[m',m].
\end{equation}
Moreover,
\begin{equation}\label{eq:107}
-\nu_{m'-1}\beta_{m'-1}+\del_{m'-1}-\del_{m'}<0
\end{equation}
(the statement being void in case that $m'=1$).
Set $\bDel_m=-\beta_m\sum_{i=m'}^m \e_i$, $\bnu'=\bnu+
\eps \e_{m}$, and $\bdel'=\bdel+\eps\bDel_m$.
The perturbation $\eps\bDel_m$ is chosen to cancel the change in $-\nu_m\beta_m$,
and to preserve the left hand side of \eqref{eq:104}
for all $i>m'-1$. However, by \eqref{eq:107}
the inequality is maintained for $i=m'-1$ as well, by taking
$|\eps|$ sufficiently small. Similarly, by \eqref{eq:108},
\eqref{eq:102} and \eqref{eq:105} are also maintained on taking
$|\eps|$ small.

This completes the construction of $\bDel_m$.
As described above, this leads to a contradiction, and we conclude that
any extreme point $(\bnu,\bdel)\in E$ satisfies $\bnu=1_k$ for some $k\in A$.

When $\bnu=1_k$ for some $k\in A$, inequalities \eqref{eq:104} can be rewritten as
$$
\del_i\le\del_{i+1},\ i\ne k,\quad \del_k\le \beta_k+\del_{k+1}.
\eqno{(\ref{eq:104}')}
$$
In particular, $0\le\del_1\le\cdots\le\del_k$.
Let $s$ denote the largest $j\le k$ for which $j\in B\cup O$
(and $s=0$ if there is no such $j$). 
The definitions of $O$ and $B$ imply $\delta_s=0$,
and thus
\stepcounter{equation}
$$
\text{$\del_j=0$ for $1\le j\le s$.}
\eqno{(\theequation a)}
$$
Similarly, $\del_{k+1}\le\cdots\le\del_J\le0$,
and if $t$ denotes the least $j\in [k+1,J]$ for which
$j\in I\cup O$ (and $t=J+1$ if empty), then
$$
\del_j=0,\quad t\le j\le J.
\eqno{(\theequation b)}
$$
Hence $\bdel$ must satisfy
$$
\begin{cases}
0\le\del_{s+1}\le\cdots\le\del_k\le\beta_k+\del_{k+1},
&
\text{when $s\le k-1$},\\
0\le\del_k
&
\text{when $s=k$}.
\end{cases}
\eqno{(\theequation c)}
$$
In addition,
$$
\begin{cases}
\del_{k+1}\le\cdots\le\del_{t-1}\le0,&
\text{when }t\ge k+2\\
\del_{k+1}\le0 &
\text{when }t=k+1.
\end{cases}
\eqno{(\theequation d)}
$$
We have just shown that \eqref{eq:102}, \eqref{eq:105}
and (\ref{eq:104}') imply (36$a$--$d$). On the other hand, clearly
(36$a$--$d$) implies (\ref{eq:104}').
Moreover, as follows directly from the definition of $s$ and $t$,
\begin{equation}\label{eq:106}
\text{$i\in I$ for all $i\in[s+1,k]$ and $i\in B$
for all $i\in[k+1,t-1]$}.
\end{equation}
This shows that (36$a$--$d$) implies \eqref{eq:102} and \eqref{eq:105}.
Thus \eqref{eq:102}, \eqref{eq:105} and (\ref{eq:104}')
are equivalent to (36$a$--$d$).
Now, the set of $\bdel$ satisfying
the constraints (36$a$--$d$) is easy to analyze.
In particular, it is not hard to see that it has the following
$t-s$ extreme points, indexed by $k$ and $r\in[s,\ldots,t-1]$, namely
\begin{equation}\label{eq:110}
\del_i^{(r,k)}=\beta_k(1_{i\ge r+1}-1_{i\ge k+1}),\quad
i=1,\ldots,J,\ r=s,\ldots,t-1.
\end{equation}
\qed

We now calculate $h$ (cf.\ \eqref{eq:109}) at each extreme point.
To this end, note that if $(\bnu,\bdel)=
(1_k,\bdel^{(r,k)})$ then by \eqref{eq:110}
$$
-\nu_i\beta_i+\del_i-\del_{i+1}=
\begin{cases}
-\beta_k & i=r,\\ 0 & i\ne r,
\end{cases}
\qquad i=1,\ldots,J,
$$
and
$$
\del_1=
\begin{cases}
\beta_k & r=0,\\ 0 & r>0.
\end{cases}
$$
Substituting in \eqref{eq:109}, we have the following possibilities.
If $r=0$ then $\del_i^{(r,k)}=\beta_k1_{i\in[1,k]}$
thus $\sum\nu_i\beta_i-\del_1=0$, and
$-\nu_i\beta_i+\del_i-\del_{i+1}=0$, $i=1,\ldots,J$. Hence
$h(1_k,\bdel^{(r,k)})=c>0$.
Otherwise,
\begin{equation}\label{eq:111}
h(1_k,\bdel^{(r,k)})=c+\la(1-e^{\beta_k})+\mu_r(1-e^{-\beta_k}).
\end{equation}
In case that $r=k$, the right hand side of \eqref{eq:111} vanishes owing to the definition of $\beta_k$ (see \eqref{eq:100}).
In case that $r<k$, recall that $k\in A$, and in particular,
$k\in A'$. By \eqref{eq:112} we therefore have $\mu_k\le\mu_r$,
hence by \eqref{eq:100}, $h(1_k,\bdel^{(r,k)})\ge0$.
Finally, consider the case where $k<r\le t-1$. By \eqref{eq:106},
$[k+1,r]\subset B$. Since $k\in A$, \eqref{eq:113} implies
that $\mu_k<\mu_r$, and we again conclude that $h(1_k,\bdel^{(r,k)})>0$.
Having shown that $h(\bnu,\bdel)\ge0$ for all
extreme points of the set $S$, we conclude that the inequality holds
for all $(\bnu,\bdel)\in S$.

\noi\uu{\it Step 2.}
We now relax the condition \eqref{eq:104}.
Let then $(\bnu,\bdel)$ satisfy \eqref{eq:102}, \eqref{eq:103},
\eqref{eq:105}, and assume that the inequality in \eqref{eq:104} holds
for all $i\in I$ (relaxing the assumption made in Step 1, that it
holds for all $1\le i\le J$). Let $P=P(\bnu,\bdel)$ denote the set of
$i$ such that
\begin{equation}\label{eq:130}
-\nu_i\beta_i+\del_i-\del_{i+1}>0.
\end{equation}
If $P$ is empty then the results of Step 1 apply and $H(\p )\ge0$.
Hence assume $P$ is not empty.
Let $j=j(\bnu,\bdel)$ be the least element in $P$.

Note that $J\not\in P$. For if $J\in P$ then using
\eqref{eq:130} and \eqref{eq:105} we find $\del_J>0$, and therefore $J\in I$. 
However, we get to assume \eqref{eq:104} for all $i \in I$, which means that 
$-\nu_J\beta_J+\del_J-\del_{J+1}\le0$.  This contradicts $J\in P$.

We proceed by backward induction on the value of $j(\bnu,\bdel)$.
Note that the sets $B$ and $I$ depend on $x$. The argument below treats
simultaneously all $x\in G$ by considering all
possible sets $B$ and $I$.

\noi
{\it Induction step.} Assumption: For all $(I,B,\bnu,\bdel)$
satisfying \eqref{eq:102}, \eqref{eq:103}, \eqref{eq:105},
such that $P(\bnu,\bdel)\cap I=\emptyset$ and $j(\bnu,\bdel)\in
[i+1,i+2,\ldots,J-1]$, one has $H(\p )\ge0$.
Let $(I,B,\bnu,\bdel)$ be such that $P(\bnu,\bdel)\cap I=\emptyset$
and $j(\bnu,\bdel)=i$ ($i\ge1$). Then
$-\nu_i\beta_i+\del_i-\del_{i+1}>0$.
Modify $\del_{i+1}$ by increasing it so as to get equality i.e.,
set $\del_{i+1}'>\del_{i+1}$ so that
$-\nu_i\beta_i+\del_i-\del'_{i+1}=0$.
This modification does not change the value of
$0\vee\mu_i(1-e^{-\nu_i\beta_i+\del_i-\del_{i+1}})$
(but keeps it zero), and it can only decrease (or leave unchanged)
the value of
$0\vee\mu_{i+1}(1-e^{-\nu_{i+1}\beta_{i+1}+\del_{i+1}-\del_{i+2}})$.
Hence by \eqref{eq:114}, the value of $H(\p )$ is only decreased.
At the same time, $j(\bnu,\bdel')>j(\bnu,\bdel)$.
Note that in modifying $\bdel$, \eqref{eq:102} need not hold
for $B$ and $I$.
However, clearly there are other sets, $B'$ and $I'$
with which it holds. For example,
\begin{equation}\label{eq:151}
B'=\{i\in[1,J]:\del'_i<0\},\quad
I'=\{i\in[1,J]:\del'_i\ge0\},\quad
O'=\emptyset.
\end{equation}
Moreover, $1\not\in B'$, since we had
$1\not\in B$ and $\del_1$ was not modified.
By the induction assumption we therefore obtain that $H(\p )\ge0$.

\noi {\it Induction base.} We show that for all $(I,B,\bnu,\bdel)$
satisfying \eqref{eq:102}, \eqref{eq:103}, \eqref{eq:105},
such that $P(\bnu,\bdel)\cap I=\emptyset$ and $j(\bnu,\bdel)=J-1$,
one has $H(\p )\ge0$.
We have $-\nu_{J-1}\beta_{J-1}+\del_{J-1}-\del_J>0$.
Similar to before, we set
$\del_J'>\del_J$ so that $-\nu_{J-1}\beta_{J-1}+\del_{J-1}-\del'_J=0$.
We claim that $P(\bnu,\bdel')=\emptyset$. Indeed, since
$P(\bnu,\bdel)=\{J-1\}$, clearly $P(\bnu,\bdel')\subset\{J\}$.
However,
\begin{align*}
-\nu_J\beta_J+\del'_J &=
-\nu_J\beta_J-\nu_{J-1}\beta_{J-1}+\del_{J-1}\\
&\le
0,
\end{align*}
where we have used the fact that $J-1\in P(\bnu,\bdel)$
implies $J-1\not\in I$, and hence by \eqref{eq:102} $\del_{J-1}\le0$.
This shows that $P(\bnu,\bdel')=\emptyset$.
As in the previous paragraph, because of the modification of $\bdel$,
\eqref{eq:102} need not hold for the $I$ and $B$ we started with,
but there are other sets $I'$ and $B'$ (defined e.g.\ as in \eqref{eq:151})
with which it holds.
The results of Step 1 therefore apply, and therefore $H(\p )\ge0$.

This completes the argument by induction and establishes
\eqref{eq:up} for superdifferentials.

\noi\uu{\it Verifying the PDE for subdifferentials.}
We are required to show that \eqref{eq:down} holds for
$\p \in D^-V(\x )$, where $\x \in G\setminus\pl_c G$ (and hence $B(\x )=\emptyset$).
Unless $D^-V(\x )$ is empty, the general form of $\p \in D^-V(\x )$, $\x \in G$ is
$$
\p =-b_k+\bdel,
$$
where $k\in A$, and
\begin{equation}\label{eq:115}
\del_i\le0,\ i\in I,\quad
\del_i=0,\ i\not\in I.
\end{equation}
Here we have used the fact that $V$ is the minimum of smooth functions,
and therefore away from the boundary $D^-V(\x )$ can have at most one element,
equal to the gradient of any minimizing function.
It suffices to show that $H(\p )\le0$ whenever $\p \cdot\bgamma_i
=-1_{i=k}\beta_k+\del_i-\del_{i+1}>0$ for all
$i\in I$. 
Using \eqref{eq:114}, \eqref{eq:115} and \eqref{eq:100},
we find
\begin{align*}
H(\p )
&=
c+\la(1-e^{\beta_k-\del_1})+\sum_{i=1}^J0\vee\mu_i(1-e^{-1_{i=k}\beta_k
+\del_i-\del_{i+1}})\\
&=
c+\la(1-e^{\beta_k-\del_1})+\sum_{i\not\in I}\mu_i(1-e^{-1_{i=k}\beta_k
-\del_{i+1}})\\
&\le
c+\la(1-e^{\beta_k})+\sum_{i\not\in I}\mu_i(1-e^{-1_{i=k}\beta_k})\\
&\le
c+\la(1-e^{\beta_k})+\mu_k(1-e^{-\beta_k})\\
&=0.
\end{align*}
\qed

\newpage
\bibliographystyle{plain}

\end{document}